\documentclass[
twoside,english]{amsart}
\PassOptionsToPackage{pdfauthor={Vladimir Kisil},%
    pdftitle={Spectrum as the Support of Functional Calculus},%
    pdfsubject={Functional Analysis},%
    backref=page,%
  pdfkeywords={symmetries, functioanl calculsu, spectrum, intertwining
  operators, jet bundle}
}{hyperref}

\let\oldtheenumi=\theenumi
\renewcommand{\theenumi}{\textup{\oldtheenumi}}
\usepackage[extnum,printedin]{myamsart}

\newcommand{\uir}[1]{\rho_{#1}}

\ppnum{LEEDS--MATH--PURE--2003-31}{arXiv: \href{http://arXiv.org/abs/math/0311285}{math/0311285}}
{2003}

\ppnum{LEEDS--MATH--PURE--2003-31
}{arXiv: \href{http://arXiv.org/abs/math/0311285}{math/0311285}}
{2003}

\input{mydef}
\providecommand{\matr}[4]{{\ensuremath{ \left(\!\! \begin{array}{cc}
#1 & #2 \\ #3 & #4
\end{array}\!\!\right) }}}

\providecommand{\Moeb}[2][\comment]{{\ensuremath{\mathcal{M}(#1,#2)}}}
\providecommand{\MoebS}[1][n]{\Moeb{\Space{B}{#1}}}

\providecommand{\oper}[1]{\mathcal{#1}}
\providecommand{\ga}{\vecbf{A}}

\providecommand{\vecbf}[1]{\mathbf{#1}}
\providecommand{\spZ}{\Space[+]{Z}{n}}

\renewcommand{\spec}{\mathrm{sp}\,}
\renewcommand{\Cliff}[2][\comment]{{\ensuremath{%
\mathcal{C}\kern-0.12em\ell(#1,#2)}}}
\ifundefined{eoe}
   \DeclareMathSymbol{\eoe}{\mathord}{AMSa}{"06}
\fi

\newcounter{myenumi}
\renewcommand{\themyenumi}{(\alph{myenumi})}

\providecommand{\SL}{\ensuremath{SL(2,\Space{R}{})}}

\usepackage{graphicx}

\begin{document}

\author[Vladimir V. Kisil]%
{\href{http://maths.leeds.ac.uk/~kisilv/}{Vladimir V. Kisil}}
\thanks{On  leave from Odessa University.}

\address{%
School of Mathematics\\
University of Leeds\\
Leeds LS2\,9JT\\
UK
}

\email{\href{mailto:kisilv@maths.leeds.ac.uk}{kisilv@maths.leeds.ac.uk}}

\urladdr{\href{http://maths.leeds.ac.uk/~kisilv/}%
{http://maths.leeds.ac.uk/\~{}kisilv/}}

\title{Monogenic
  Calculus as an Intertwining Operator}

\begin{abstract}
  We revise a monogenic calculus for several non-commuting operators,
  which is defined through group representations.  Instead of an
  algebraic homomorphism we use group covariance.
  The related notion of joint
  spectrum and spectral mapping theorem are discussed. The construction is
  illustrated by a simple example of calculus and joint spectrum of two
  non-commuting selfadjoint \(n\times n\) matrices.
\end{abstract}
  \keywords{Functional calculus, spectrum, intertwining operator,
    spectral mapping theorem, jet spaces, monogenic function, Clifford
  algebra}
  \AMSMSC{47A60}{30G35, 46H30, 47A10, 47B15}
\maketitle

\tableofcontents\vspace{3mm}

\section{Introduction}
\label{sec:introduction-1}

Central objects of operator theory are functional calculus (usually
defined as an algebra homomorphism), spectrum (defined as set of
singular points of the resolvent), and spectral mapping theorem
(describing transformations of the spectrum under the functional
calculus). Following the discussion in~\cite{Kisil02a} we arrange
these objects as follows: 
\begin{enumerate}
\item Functional  calculus is an \emph{original} notion defined in
  some independent terms;
\item Spectrum is derived from the previously
  defined functional calculus as its \emph{support} in some
  appropriate sense;
\item Then Spectral mapping theorem should drop out naturally.
\end{enumerate}

The full potential of such a construction depends from its
source---the definition of a functional calculus. It is known that
homomorphic calculi are successful only in few simplest cases,
e.g. for a single normal operator. To increase its potential
functional calculus through group covariance was defined
in~\cite{Kisil95i}. In that paper such a definition was applied to a
complicated case of monogenic calculus for several non-commuting
operators. Monogenic calculus of commuting operators was considered
earlier in~\cite{McInPryde87}. Monogenic calculus of non-commuting
operators was carefully developed through plane wave decomposition in
many subsequent papers,
see~\cite{JeffMcInt98a,JefferiesMcIntosh99a,JefferiesJohnson01,JefferiesStraub03}. These
papers contain many important results, e.g. connection between
monogenic and Weyl~\cite{Anderson69} calculi, but they do not consider
covariant properties of the
calculus. Papers~\cite{KisRam95a,KisRam96a} utilised an algebraic
approach and do not develop group representations either.

Meanwhile the covariant approach to functional calculus required
development of wavelets technique and its applications to analytic
function theory, it was performed
in~\cite{Kisil97c,Kisil97a,Kisil98a,Kisil01a,Kisil02c}.  It emerged
from these studies that the new definition is a useful replacement for
classical one across all range of problems, even in case of a single
non-normal operator with finite range~\cite{Kisil02a}. The key
ingredient in this approach is the development of all principal
objects of analytical function theory (Cauchy integral, Hardy and
Bergman spaces, Cauchy-Riemann equations, Taylor series, etc.) from
the group of M\"obius transformations and wavelets
technique~\cite{Kisil97c,Kisil01a}. This allows to give a template
definition of functional calculus as follows,
cf.~\cite[Defn.~1.1]{Kisil95i}:  

\begin{defn}
  \label{de:functional-calculus-new}
  Let \(A\) be a normed algebra,  and  \(M\)\ be a left
  \(\algebra{A}\)-module. Let \(G\) be a group, \(X\) be a left
  \(G\)-homogeneous space, and \(\FSpace{A}{}(X)\) be an associated
  space of analytic functions.
  An \emph{analytic functional calculus} for an element
  \(a\in\algebra{A}\)\  is a \textit{continuous 
  linear} mapping
  \(\Phi:\FSpace{A}{}(X)\rightarrow \FSpace{A}{}(X,M)\)\ such that 
  \begin{enumerate} 
  \item \(\Phi\) is an \emph{intertwining operator} 
    \begin{displaymath}
      \Phi {\uir{X}} = {\uir{M}} \Phi
    \end{displaymath}
    between two representations of the
    group \(G\): \(\uir{X}\) acts in the analytic space
    \(\FSpace{A}{}(X)\) of scalar valued functions on \(X\) and \(\uir{M}\)\
    acts in a space \(\FSpace{A}{}(X,M)\) of \(M\)-valued functions in
    a way depending from \(a\in\algebra{A}\).
  \item There is an initialisation condition: \(\Phi[f_0]=f_M\),
    i.e. the vacuum vector of \(\FSpace{A}{}(X)\) is mapped into the
    vacuum vector of \(\FSpace{A}{}(X,M)\).  
  \end{enumerate}
\end{defn} Note that our functional calculus released from the
homomorphism condition can take value in any left
\(\algebra{A}\)-module \(M\), which however could be \(\algebra{A}\)
itself if suitable. This improves spectral localisation technique in
our construction.

In the paper~\cite{Anderson69}
joint spectrum was defined as the support of the Weyl calculus, i.e. as
the set of points where the operator valued distribution does not
vanish. We also define
the spectrum as a support of functional calculus, but due to our
Definition~\ref{de:functional-calculus-new} it has a different
meaning. 
\begin{defn}\textup{\cite{Kisil02a}}
  \label{de:spectrum-new}
    A corresponding \emph{spectrum} of \(a\in\algebra{A}\) is the
  \textit{support} of the functional calculus \(\Phi\), i.e. the
  collection of non-vanishing intertwining operators between
  \(\uir{M}\) and \emph{prime
  representations}~\cite[\S~8.3]{Kirillov76}.
\end{defn}

More variations of functional calculi (Weyl, Wick, Berezin, etc.) are
obtained from other groups and their
representations~\cite{Kisil95i,Kisil98a}. There are also recent papers
of other researchers devoted to covariant
calculus~\cite{ArazyUpmeier02a,ArazyUpmeier02b}.

\section{Preliminaries on Clifford Algebras and M\"obius Transformations }
\label{sec:cliff-algebr-mobi}

Let \(\Space{R}{n}\) be a real \(n\)-dimensional vector space with a
fixed frame \(e_1\), \(e_2\), \ldots, \(e_n\). 
Let \(\Cliff{n}\) be the \emph{real Clifford
  algebra} generated by \(1\), \(\vecbf{e}_j\), \(1\leq
j\leq n\)
and the relations 
\begin{displaymath}
  \vecbf{e}_i \vecbf{e}_j + \vecbf{e}_j \vecbf{e}_i =-2\delta_{ij}.
\end{displaymath} Then there is the natural embedding of
\(\Space{R}{n}\) into \(\Cliff{n}\). We identify \(\Space{R}{n}\) with
its image in \(\Cliff{n}\) and call its elements
\emph{vectors}. There are two linear anti{-}automorphisms \(*\)
(reversion) and \(\bar\ \) (main anti{-}automorphisms) and automorphism
\('\) of \(\Cliff{n}\) defined on its basis
\(A_\nu=\vecbf{e}_{j_1}\vecbf{e}_{j_2}\cdots \vecbf{e}_{j_r}\),
\(1\leq j_1 <\cdots<j_r\leq n\) by the rule:
\begin{eqnarray*}
  (A_\nu)^*= (-1)^{\frac{r(r-1)}{2}} A_\nu, \qquad
  \bar{A}_\nu= (-1)^{\frac{r(r+1)}{2}} A_\nu,\qquad
  A_\nu'= (-1)^{r} A_\nu.
\end{eqnarray*}
In particular, for vectors, \(\bar{\vecbf{x}}=\vecbf{x}'=-\vecbf{x}\) and
\(\vecbf{x}^*=\vecbf{x}\).

It is easy to see that \(\vecbf{x}\vecbf{y}=\vecbf{y}\vecbf{x}=1\)
for any \(\vecbf{x}\in\Space{R}{n}\) such that
\(B(\vecbf{x},\vecbf{x})\neq 0\)
and \(\vecbf{y}={\bar{\vecbf{x}}}\,{\norm{\vecbf{x}}^{-2}}\), which is
the \emph{Kelvin inverse} of \(\vecbf{x}\).
Finite products of invertible vectors are invertible in \(\Cliff{n}\)
and form the \emph{Clifford group} \(\Gamma(n)\)~\cite[(1.39)]{Cnops02a}. Elements
\(a\in\Gamma(n)\) such that
\(a\bar{a}=\pm 1\) form the \(\object[(n)]{Pin}\) group---the double cover
of the group of orthogonal rotations \(\object[(n)]{O}\). We also
consider~\cite[\S~5.2]{Cnops94a} \(T(n)\) to be the set of all
products of vectors in \(\Space{R}{n} \).

Let \((a, b, c, d)\) be a quadruple from \(T(n)\) with
the properties:
\begin{enumerate}
\item \((ad^*-bc^*)\in \Space{R}{}\setminus {0}\);
\item \(a^*b\), \(c^*d\), \(ac^* \), \( bd^*\) are vectors.
\end{enumerate}
Then~\cite[Thm.~5.2.3]{Cnops94a}, \cite[(4.10)]{Cnops02a}
\(2\times 2\)-matrices \matr{a}{b}{c}{d} form the
group \(\Gamma(1,n+1)\) under the
usual matrix multiplication. It has a representation
\(\uir{\Space{R}{n} }\) \comment{we denote its restriction to any
  subgroup by the same notation.} by transformations of
\(\dot{\Space{R}{n} }\) given by:
\begin{equation}\label{eq:sp-rep}
  \uir{\Space{R}{n}}\matr{a}{b}{c}{d} :
  \vecbf{x} \mapsto (a\vecbf{x}+b)(c\vecbf{x}+d)^{-1},
\end{equation}
which form the \emph{M\"obius} (or
the \emph{conformal}) group of \(\dot{\Space{R}{n}}\). 
Here \(\dot{\Space{R}{n} 
}\) the compactification of \(\Space{R}{n} \) by the point  at infinity 
(see~\cite[\S~5.1]{Cnops94a}).
The analogy with fractional-linear transformations of the complex line
\Space{C}{} is useful, as well as representations of shifts
\(\vecbf{x}\mapsto \vecbf{x}+y\), orthogonal rotations
\(\vecbf{x}\mapsto k(a)\vecbf{x}\), dilations
\(\vecbf{x}\mapsto \lambda \vecbf{x}\), and the Kelvin inverse
\(\vecbf{x}\mapsto \vecbf{x}^{-1}\) by the
matrices \matr{1}{y}{0}{1}, \matr{a}{0}{0}{{a}^{*-1}},
\matr{\lambda^{1/2}}{0}{0}{\lambda^{-1/2}}, \matr{0}{-1}{1}{0}
respectively.

Following e.g.~\cite[(1.8)]{Cnops02a} we adopt the next agreement.
\begin{notation}
  \label{no:division}
  In a non-commutative algebra setting the ambiguous notation
  \(\displaystyle\frac{a}{b}\) always means \(ab^{-1}\). Consequently 
  \(\displaystyle\frac{ac}{bc}=\displaystyle\frac{a}{b}\) but
  \(\displaystyle\frac{ca}{cb}\neq \displaystyle\frac{a}{b}\) in
  general.
\end{notation}

Study of M\"obius transformation is facilitated by introduction of
projective coordinates in the space \(P\Space{R}{1,n+1}\) of spheres
in \(\Space{R}{n}\) \cite[(4.12)]{Cnops02a}.  The spere with the
centre \(m\in\Space{R}{n}\) and the radius \(r\) defined by the
equation \((\vecbf{y}-\vecbf{m})^2=r^2\) is associated with the ray of
matrices by the map
\begin{equation}
  \label{eq:spheres-to-matrices}
T: \quad  \{ \vecbf{y} \such B(\vecbf{y}-\vecbf{m},\vecbf{y}-\vecbf{m})=r^2 \} 
 \quad\mapsto \quad
\lambda \begin{pmatrix}
  \vecbf{m} & -\vecbf{m}^2-r^2\\ 1 & -\vecbf{m}
\end{pmatrix}.
\end{equation}
 A point \(\vecbf{x}\in\Space{R}{n}\) is associated with a
zero radius sphere with the centre \(\vecbf{x}\) and thus is represented by \(
\begin{pmatrix}
  \vecbf{x} & -\vecbf{x}^2\\ 1 & \vecbf{x}
\end{pmatrix}\). Then M\"obius transformations~\eqref{eq:sp-rep}
corresponds to the orthogonal rotations in the projective space
\(P\Space{R}{1,n+1}\) as follows \cite[(4.13)]{Cnops02a}:
\begin{equation}
  \label{eq:project-ortogonal}
  \uir{P} \begin{pmatrix}
    a&b\\c&d
  \end{pmatrix}: \quad
\begin{pmatrix}
  \vecbf{m} & -\vecbf{m}^2-r^2\\ 1 & -\vecbf{m}
\end{pmatrix} \mapsto
  \begin{pmatrix}
    a&b\\c&d
  \end{pmatrix}
\begin{pmatrix}
  \vecbf{m} & -\vecbf{m}^2-r^2\\ 1 & -\vecbf{m}
\end{pmatrix}
  \begin{pmatrix}
    \bar{d}&\bar{b}\\\bar{c}&\bar{a}
  \end{pmatrix}.
\end{equation}

One usually says that the conformal group in \(\Space{R}{n}\), \(n>2\) is
not so rich as the conformal group in \(\Space{R}{2}\).  Nevertheless,
the conformal covariance has many applications in Clifford
analysis~\cite{Cnops94a,Ryan95b}.  Notably, groups of conformal
mappings of unit spheres \(\Space{S}{n-1}=\{\vecbf{x} \such
\vecbf{x}\in\Space{R}{n}, B(\vecbf{x},\vecbf{x})=1 \} \) onto itself
are similar for all \(n\) and as sets can be parametrised by the
product of the unit ball \(\Space{B}{n}\) and
the group of isometries of \Space{S}{n-1}. We specialise main result of
\cite{Kisil96d} for the positive definite case as follow:
\begin{prop} 
  The group \(\MoebS[n]\) of conformal mappings of the open unit
  sphere \(\Space{S}{n-1}\) onto itself represented by matrices 
  \begin{equation}
    \label{eq:matrices-alpha-beta}
    \matr{a }{b'}{b}{a '}, \qquad a ,b\in 
    T(n),\quad ab^*\in \Space{R}{n},\quad  \modulus{a}^2-\modulus{b}^2=1.
  \end{equation}
  Its inverse is \(\matr{a }{b'}{b}{a '}^{-1}=\matr{\bar{a}}{-\bar{b}}{-b^*}{a^*}\).
\end{prop}
\begin{proof}
  The proof is easy in the projective
  coordinates~\eqref{eq:spheres-to-matrices}. Indeed the unit spere
  corresponds to the matrix \(
  \begin{pmatrix}
    0&-1\\1&0
  \end{pmatrix}\). Straightforwardly transformations
  \(\uir{P}\)~\eqref{eq:project-ortogonal} with matrices of the
  form~\eqref{eq:matrices-alpha-beta} preserve this ray:
  \begin{displaymath}
    \matr{a}{b'}{b}{a'}\matr{0}{-1}{1}{0}\matr{a^*}{b^*}{\bar{b}}{\bar{a}} =
    \matr{b'}{-a}{a'}{-b}\matr{a^*}{b^*}{\bar{b}}{\bar{a}}=\matr{0}{-1}{1}{0}.
  \end{displaymath}
 Thus corresponding M\"obius transformations preserve the unit sphere.
\end{proof}

The presentation of \(\MoebS\) by~\eqref{eq:matrices-alpha-beta} is
difficult to use due to ineffective definition through the constrain
\(\modulus{a}^2-\modulus{b}^2=1\). Thus we will prefer a direct
parametrisation, cf.~\cite[\S~VI.1.3]{Vilenkin68}, as follows.
We can identify the unit ball \(\Space{B}{n}\) with the left coset
\(O(\Space{S}{n-1}) \backslash \MoebS\), where the decomposition
\(\MoebS\sim O(\Space{S}{n-1})\times\Space{B}{n}\) follows
from~\eqref{eq:u-w-form-of-Mn}. Note that \(K=O(\Space{S}{n-1})\) is the
maximal compact subgroup of \(\MoebS\).
\begin{eqnarray}
  \label{eq:sl2-u-psi-coord}
  \matr{a}{b'}{b}{a'} 
  & =& \modulus{a} 
  \matr{ 
    \frac{a}{ \modulus{a} } }{0}{0}{\frac{a'}{ 
      \modulus{a} } }
  \matr{1}{\frac{\bar{a}}{ \modulus{a}^2}b'}{\frac{a^*}{\modulus{a}^2}b}{1}
  \nonumber \\
  &=& \frac{1}{\sqrt{1+\vecbf{u}^2 }}
  \matr{w}{0}{0}{w'}
  \matr{1}{\vecbf{u}'}{\vecbf{u}}{1}, 
\end{eqnarray}
where
\begin{equation}
  \label{eq:connection-ab-to-uw}
  w=\frac{a}{\modulus{a}},\qquad 
  \vecbf{u}=\frac{a^*}{\modulus{a}^2}b,\qquad
  \sqrt{1+\vecbf{u}^2 }=\modulus{a}^{-1},\qquad 
  \modulus{\vecbf{u}}<1.
\end{equation}
Consequently for \(\vecbf{u}\in\Space{B}{n}\), \(w\in\Gamma(n)\)
the M\"obius transformations \(\phi_{(\vecbf{u},w)}\) with
matrix
\begin{equation}
  \label{eq:u-w-form-of-Mn}
  \frac{1}{\sqrt{1+\vecbf{u}^2}} \matr{w}{0}{0}{w'}\matr{1}{\vecbf{u'}}{\vecbf{u}}{1} 
  = \frac{1}{\sqrt{1+\vecbf{u}^2}}\matr{w}{w\vecbf{u}'}{w'\vecbf{u}}{w'}, 
\end{equation} constitute \(\MoebS[n]\). Sometime in M\"obius
transformations we will omit the normalising factor
\((1+\vecbf{u}^2)^{-1/2}\) in~\eqref{eq:u-w-form-of-Mn} for the sake
of brevity. Although one should not forget that this factor is
important for the invariant integration on \(\MoebS\).
The following properties follows from such a realisation: 
  \begin{lem}
    \(\MoebS[n]\) acts on \(\Space{B}{n}\) transitively.
    Transformations of the form \(\phi_{(0,w)}\) constitute a subgroup
    isomorphic to \(\object[(n)]{O}\). The homogeneous space
    \(\MoebS[n]/\object[(n)]{O}\) is isomorphic as a set to
    \(\Space{B}{n}\). Moreover:
  \begin{enumerate}
  \item\label{it:ident1} \(\phi_{(\vecbf{u},1)}^2=-1\) on
    \(\Space{B}{n}\), thus
    \(\phi_{(\vecbf{u},1)}^{-1}=\phi_{(\vecbf{u}',1)}=\phi_{(-\vecbf{u},1)}\).\vspace{3pt}
  \item \(\phi_{(\vecbf{u},1)}^{-1}(0)=-\vecbf{u}\) and
    \(\phi_{(\vecbf{u},1)}^{-1}(\vecbf{u}\strut)=0\).\vspace{3pt}
  \item\label{it:ident3}
    \(\phi_{(\vecbf{u}_1,1)}^{-1}\phi_{(\vecbf{u}_2,1)}^{-1}=
    \phi_{(\vecbf{u},w)}^{-1}\) \vspace{3pt} where
    \begin{displaymath}
    \vecbf{u}=\phi_{(\vecbf{u}_1,1)}^{-1}(\vecbf{u}_2)
    =\phi_{(\vecbf{u}_2,1)}^{-1}(\vecbf{u}_1)\quad \textrm{ and } \quad
    w=\frac{1-\vecbf{u}_1\vecbf{u}_2}{\modulus{1-\vecbf{u}_1\vecbf{u}_2}}.
  \end{displaymath}
  \end{enumerate}
\end{lem} We use the same notation for the M\"obius transformation
\(\phi_{(\vecbf{u},w)}\) and the matrix~\eqref{eq:u-w-form-of-Mn}
which produces it. It is a direct check to see that 
\begin{displaymath}
  \matr{w}{0}{0}{w'}\matr{1}{\vecbf{u}'}{\vecbf{u}}{1}  
=   \matr{1}{w\vecbf{u}'w^*}{w'\vecbf{u}\bar{w}}{1}\matr{w}{0}{0}{w'},
\end{displaymath}
which implies that
\(
  \phi_{(\vecbf{u},w)}^{-1}=\phi_{({w}^*\vecbf{u}'w,\bar{w})}.
\)

\begin{lem}
  \label{le:haar}
  The left invariant \emph{Haar measure} \(dg\) 
  on \(\MoebS\sim \Space{B}{n}\times O(n)\) in coordinates \((\vecbf{u},w)\)  is
  \begin{equation}
    \label{eq:Haar-measure}
    dg(\vecbf{u},w)=\frac{d\vecbf{u}\,dw}{\modulus{1+\vecbf{u}^2}^{n}},
  \end{equation}
  where \(dw\) is a Haar measure on
  \(O(n)\) and \(d\vecbf{u}\) is the Lebesgue measure on
  \(\Space{B}{n}\).
\end{lem}
\begin{proof}
  It follows from~\ref{it:ident3} that left shifts on \(\MoebS\) in
  coordinates \((\vecbf{u},w)\) acts by M\"obius transformations on
  \(\vecbf{u}\in\Space{B}{n}\) which fixes the unit sphere. According
  to~\cite[Cor.~6.1.2]{Cnops94a} the invariant metric on
  \(\Space{B}{n}\) is defined through the distance of the zero radius
  sphere defined by \(\vecbf{u}\) to the unit sphere
  \(\Space{S}{n-1}\). This distance is
  \(\modulus{1+\vecbf{u}^2}^{-1}\), thus the invariant measure is
  obtained from its \(n\)-th power.
\end{proof}

The importance of the Haar measure is justified by the \emph{invariant
integration} (or \emph{invariant functional}) it produces:
\begin{displaymath}
  \int_{\MoebS}f(g)\,dg=\int_{\MoebS}f(g_1g)\,dg, \quad \textrm{ for
    all } f(g)\in\Space{L}{1}(\MoebS) \textrm{ and } g\in\MoebS.
\end{displaymath}
It is rarely realised that the Haar invariant functional is not
the only possible and that other invariant functionals are useful as
well. The classic example is described below.
\begin{lem}
  \label{le:hardy-functional}
  The invariant functional \(H\) on \(\MoebS\) of \emph{Hardy type} is
  given by: 
  \begin{equation}
    \label{eq:hardy-functional}
    H(f)=\lim_{r\rightarrow
    1}\int_{O(n)}\int_{\Space{S}{n-1}} f(r\vecbf{u},w)
    \,\frac{dw\,d\vecbf{u}}{\modulus{1+\vecbf{u}^2}^{n-1}},
    \quad\textrm{ where } w\in O(n), \vecbf{u}\in\Space{S}{n-1}.
  \end{equation}
\end{lem}
\begin{proof}
  This result follows from the discussion in the proof of
  Lemma~\ref{le:haar} and observation that the limit
  in~\eqref{eq:hardy-functional} is  M\"obius
  invariant. 
\end{proof}
\begin{defn}
  The \emph{Hardy inner product} in a space of Clifford valued
  functions on \(\MoebS\) is  derived from the Hardy
  functional~\eqref{eq:hardy-functional}:
  \begin{equation}
    \label{eq:Hardy-inner-prod}
    \scalar{f_1}{f_2} = H(\bar{f}_1f_2).
  \end{equation}
\end{defn}
Note that the invariance of the Hardy
functional~\eqref{eq:hardy-functional} implies that left
shifts are isometries with respect to norm defined
through~\eqref{eq:Hardy-inner-prod}. 

\section{Construction of Clifford Analysis from $\MoebS$ Group}
\label{sec:backgr-compl-analys}

\subsection{Wavelet Transform and Cauchy Kernel}
\label{sec:wavel-transcf-cauchy}

To understand the functional calculus from
Definition~\ref{de:functional-calculus-new} we need first to realise
the function theory of monogenic functions from the representation
theory of \(\MoebS\), see~\cite{Kisil97c,Kisil97a,Kisil01a,Kisil02c}
for more details. 

Each element \(g\in\MoebS\) acts by the linear-fractional transformation
(the M\"obius map) on \(\Space{B}{n}\)
and \(\Space{S}{n-1}\) \emph{from the left} as follows: 
\begin{equation}
  \label{eq:moebius}
  g^{-1}: \vecbf{x} \mapsto \frac{\bar{a} \vecbf{x} - \bar{b}}{a^*-b^* \vecbf{x}},
  \qquad \textrm{ where } \quad
  g^{-1}=\matr{\bar{a}}{-\bar{b}}{-b^*}{a^*}.
\end{equation}
In the decomposition~\eqref{eq:sl2-u-psi-coord} the first matrix on
the right hand side acts by transformation~\eqref{eq:moebius} as an
orthogonal rotation of \(\Space{S}{n-1}\) and \(\Space{B}{n}\); and the
second one---by transitive family of maps of the unit ball onto
itself.

M\"obius transformations~\eqref{eq:moebius} could be linearised to the
representation \(\uir{1}\) on functions, cf.~\cite[(4.56)]{Cnops02a}
and~\cite[Thm.~5.4.1]{GilbertMurray91}, by the induced representation
technique~\cite[\S~13]{Kirillov76}:
\begin{equation}
  \label{eq:rho-1-1}
  \uir{1}(g): f(z) \mapsto
  \frac{a'-\bar{\vecbf{x}}b'}{\modulus{a'-\bar{\vecbf{x}}b'}^n} \,
  f\!\left(
    \frac{\bar{a} \vecbf{x} - \bar{b}}{a^*- b^*\vecbf{x}} 
  \right),
\quad  \textrm{ where } \quad 
g^{-1}=\matr{\bar{a}}{-\bar{b}}{-b^*}{a^*}.
\end{equation}
Indeed one can directly verify:
\begin{eqnarray*}
  \lefteqn{\uir{1}(g_1)( \uir{1}(g_2) f(z))=  \uir{1}(g_1)\left(
      \frac{a'_2-\bar{\vecbf{x}}b_2'}{\modulus{a_2'-\bar{\vecbf{x}}b'_2}^n} \, 
      f\!\left(\frac{\bar{a}_2 \vecbf{x} - \bar{b}_2}{a_2^*- b_2^*\vecbf{x}} 
      \right)\right) } \strut\\
  &=&  \frac{a'_1-\bar{\vecbf{x}}b_1'}{\modulus{a_1'-\bar{\vecbf{x}}b'_1}^n}
  \frac{a'_2-\overline{(\bar{a}_1 \vecbf{x} - \bar{b}_1)(a_1^*-b_1^*
      \vecbf{x})^{-1}}b_2'} 
  {\modulus{a_2'-\overline{(\bar{a}_1 \vecbf{x} - \bar{b}_1)
        (a_1^*-b_1^* \vecbf{x})^{-1}}b'_2}^n}  \strut \\
  &&\quad{}\times
  f\!\left(\frac{\bar{a}_2(\bar{a}_1 \vecbf{x} - \bar{b}_1)(a_1^*-b_1^* \vecbf{x})^{-1}  -
      \bar{b}_2}
    {(a_2^*(a_1^*-b_1^* \vecbf{x})(a_1^*-b_1^* \vecbf{x})^{-1}- b_2^*}  
  \right)  \strut\\
  &=&  
  \frac{(a'_1-\bar{\vecbf{x}}b_1')a'_2-( \bar{\vecbf{x}}{a}_1 - {b}_1)b_2'} 
  {\modulus{(a'_1-\bar{\vecbf{x}}b_1')a'_2-( \bar{\vecbf{x}}{a}_1 - {b}_1)b_2'}^n} \,
  f\!\left(\frac{\bar{a}_2(\bar{a}_1 \vecbf{x} - \bar{b}_1)  -
      \bar{b}_2(a_1^*-b_1^* \vecbf{x})}
    {a_2^*(a_1^*-b_1^* \vecbf{x})- b_2^*(\bar{a}_1 \vecbf{x} - \bar{b}_1)} 
  \right) \strut\\
  &=&  
  \frac{(a'_1a'_2+{b}_1b_2')-\bar{\vecbf{x}}(b_1'a'_2+{a}_1b_2') } 
  {\modulus{(a'_1a'_2+{b}_1b_2')-\bar{\vecbf{x}}(b_1'a'_2+{a}_1b_2')}^n} \,
  f\!\left(\frac{(\bar{a}_2\bar{a}_1+\bar{b}_2b_1^*) \vecbf{x} - (\bar{a}_2\bar{b}_1 +
      \bar{b}_2a_1^*)}
    {(a_2^*a_1^*+b_2^* \bar{b}_1)-(a_2^*b_1^* + b_2^*\bar{a}_1) \vecbf{x} } 
  \right) \strut \\
  &=&  
  \frac{a'-\bar{\vecbf{x}}b' } {\modulus{a'-\bar{\vecbf{x}}b'}^n} \,
  f\!\left(\frac{\bar{a}\vecbf{x} - \bar{b}}{a^*-b^* \vecbf{x}
    }\right)
  \quad  [\textrm{where}\ a=a_1a_2+b_1'b_2, b=b_1a_2+a_1'b_2] \strut\\
  &=& \uir{1}(g_1g_2) f(\vecbf{x}).  \strut
\end{eqnarray*}

Let \(\FSpace{L}{2}(\Space{S}{n-1})\) be equipped with a Clifford valued
inner product, cf.~\cite[(1.29)]{Cnops02a}:
\begin{equation}
  \label{eq:inner-product}
  \scalar{f_1}{f_2}=\int_{\Space{S}{n-1}}
  \bar{f}_1(\vecbf{x})f_2(\vecbf{x})\, d\vecbf{x}
\end{equation} normalised such that
\(\int_{\Space{S}{n-1}}d\vecbf{x}=1\).  Then~\cite[(4.56)]{Cnops02a}
the representation~\eqref{eq:rho-1-1} became unitary in
\(\FSpace{L}{2}(\Space{S}{n-1})\).

We choose~\cite{Kisil97c,Kisil98a,Kisil01a} \(K\)-invariant function
\(f_0(\vecbf{x})=(\vecbf{x})\equiv 1\) be \emph{vacuum vector} or 
\emph{mother wavelet}~\cite{Kisil98a}. Then \emph{coherent states} or
\emph{wavelets} are all transformations of the vacuum vector by
\(\uir{1}\):
\begin{equation}
  \label{eq:wavelets}
  f_g(\vecbf{x})=\uir{1}(g)f_0(\vecbf{x})=
  \frac{a'-\bar{\vecbf{x}}b'}{\modulus{a'-\bar{\vecbf{x}}b'}^n},
  \qquad g^{-1}=\matr{\bar{a}}{-\bar{b}}{-b^*}{a^*}.  
\end{equation}
They are mainly determined by the point on the unit disk \(
\vecbf{u}=a^*b/\modulus{a}^2\). 
The linear span of all wavelets is called the \emph{Hardy space}
\(\FSpace{H}{2}(\Space{S}{n-1})\),  and \(f_0\) is \emph{cyclic} in
\(\FSpace{H}{2}(\Space{S}{n-1})\). 
M\"obius transformations provide a natural family of
intertwining operators for \(\uir{1}\) coming from inner
automorphisms of \(\MoebS\) (will be used later). 
 
The \emph{wavelet transform}~\cite{Kisil97c,Kisil98a}
\(\oper{W}:\FSpace{L}{2}(\Space{S}{n-1})\rightarrow
  \FSpace{H}{2}(\MoebS)\) is defined by:
\begin{eqnarray}
  \label{eq:cauchy}
  \label{eq:wavelet-transform-abstract}
  \quad\oper{W}f(g)&=&\scalar{f_g}{f} \\
  &=&\int_{\Space{S}{n-1}}
  \frac{{a}^*-{b}^*{\vecbf{x}}}
  {\modulus{{a}^*-{b}^*{\vecbf{x}}}^n} f(\vecbf{x})\, d\vecbf{x}
  \label{eq:wavelet-transf}\\
  &=&\int_{\Space{S}{n-1}}
  \frac{{a}^*\bar{\vecbf{x}}-{b}^*}{\modulus{{a}^*\bar{\vecbf{x}}-{b}^*}^n} 
  \,{\vecbf{x}} d\vecbf{x}\,f(\vecbf{x})  .
  \nonumber  \\
  &=&\frac{{a}^*}{\modulus{a}^n}\int_{\Space{S}{n-1}}
  \frac{\bar{\vecbf{x}}-\bar{\vecbf{u}}}{\modulus{{\vecbf{x}}-\vecbf{u}}^n} 
  \,d\sigma(\vecbf{x})\,f(\vecbf{x}),  \  \textrm{ where }
  \vecbf{u} =\frac{{b}'a^*}{\modulus{a}^2}, \ 
  d\sigma(\vecbf{x})={\vecbf{x}} d\vecbf{x}  . 
  \label{eq:cauchy-formula}
\end{eqnarray}
If we consider the \emph{reduced wavelet
  transform}~\cite{Kisil97c,Kisil98a} \(\oper{W}:\FSpace{L}{2}(\Space{S}{n-1})\rightarrow
  \FSpace{H}{2}(\Space{B}{n})\) then  
the last formula is the Cauchy integral formula in Clifford
analysis up to the factor \(\frac{{a}^*}{\modulus{a}^n}\). This factor
is similar to the factor \(\sqrt{1-\modulus{u}^2}\) in the Cauchy
formula in complex analysis derived in~\cite[(3.20)]{Kisil97c}. Their
appeared due to the invariant measures on \(\SL\) and \(\MoebS\). 
Note the appearance of the important Clifford valued differential form
\(d\sigma(\vecbf{x})={\vecbf{x}} d\vecbf{x}\) in~\eqref{eq:cauchy-formula},
cf.~\cite[\S~9.1]{Deinze93}, \cite[\S~II.0.2.1]{DelSomSou92}.
A standard derivation of the Cauchy formula in Clifford
analysis are based on Stokes's Theorem.  

Although the Cauchy formula (i.e. reduced wavelet transform) is an
established tool in analytic function theory its unreduced
version~\eqref{eq:wavelet-transf} acting
\(\oper{W}:\FSpace{L}{2}(\Space{S}{n-1})\rightarrow
\FSpace{H}{2}(\MoebS)\) is also valuable for the functional calculus
of several non-commuting operators.

The wavelet transform of the vacuum vector \(f_0\)
\begin{eqnarray}
    \oper{W}f_0(g)=\scalar{f_g}{f_0} &=& \int_{\Space{S}{n-1}}
  \frac{{a}^*-{b}^*{\vecbf{x}}}
  {\modulus{{a}^*-{b}^*{\vecbf{x}}}^n} \,
  d\vecbf{x}=\frac{{a}^*}{\modulus{a}^n}, \quad \textrm{ where }
  g=\matr{a}{b'}{b}{a'}, \textrm{ or}\nonumber \\
  \label{eq:vacuum-vector}
  &=&w^*(1+\vecbf{u}^2)^{(n-1)/2}, \quad \textrm{ where }
  g= \frac{1}{\sqrt{1+\vecbf{u}^2}}\matr{w}{w\vecbf{u}'}{w'\vecbf{u}}{w'}.
\end{eqnarray}
Consequently \(\oper{W}f_0\) has a finite norm with respect
to~\eqref{eq:Hardy-inner-prod}. 
\begin{defn}
  The \emph{Hardy space} \(\FSpace{H}{2}(\MoebS)\) of Clifford valued
  functions on \(\MoebS\)  is a left \(\Cliff{n}\)-module invariant
  under left shifts on \(\MoebS\), which is generated by the
  vacuum vector \(\oper{W}f_0\) \eqref{eq:vacuum-vector}.
\end{defn}
From the general wavelet technique~\cite{Kisil98a} we obtain the
following result:
\begin{lem}
  \begin{enumerate}
  \item \(\FSpace{H}{2}(\MoebS)\) is an inner product space with
    the product derived from the Hardy
    functional~\eqref{eq:hardy-functional}:
    \begin{equation}
      \label{eq:hardy-inner-product}
      \scalar{f_1}{f_2} = H(\bar{f}_1f_2), \qquad 
      \textrm{ where }
      f_1, f_2\in\FSpace{H}{2}(\MoebS).
    \end{equation}
  \item Wavelet transform~\eqref{eq:wavelet-transf} is a unitary
    operator intertwining the
    representation \(\uir{1}\) on \(\FSpace{H}{2}(\Space{S}{n-1})\) and
    the left regular representation on \(\FSpace{H}{2}(\MoebS)\) by
    shifts:
    \begin{displaymath}
      \oper{W}\uir{1}(g)  = \lambda(g)\oper{W}, \qquad \textrm{ for
        all } g\in\MoebS.
    \end{displaymath}
  \end{enumerate}
\end{lem}

\subsection{Taylor Series}
\label{sec:taylor-series}

Other classical objects of Clifford analysis (the Cauchy-Riemann
equation, the Bergman space, etc.)  can be also
obtained~\cite{Kisil97c,Kisil01a} from representation
\(\uir{1}\). However we need only the Taylor series in the present
paper. It is known~\cite[\S~11.2.2]{BraDelSom82} that there is the
orthonormal basis \(V_m(\vecbf{x})\) of
\(\FSpace{H}{2}(\Space{S}{n-1})\) labelled by a multiindex
\(m=(m_1,\ldots,m_n) \in\spZ\). Elements \(V_m(\vecbf{x})\)
can be constructed as symmetric polynomials of hypercomplex variables
\(e_1x_j-e_jx_1\), \(j=1,\ldots,n\). Consequently there is a
decomposition of the Cauchy kernel (i.e. coherent
states~\eqref{eq:wavelets}):
\begin{equation}
  \label{eq:Cauchy-decomposition}
  \uir{1}(g)f_0(\vecbf{x})=f_g(\vecbf{x})=\sum_{m\in\spZ}
  W_m(g)\,V_m(\vecbf{x}), 
  \quad \textrm{ where } g\in\MoebS
\end{equation}
with some functions on \(\MoebS\) defined by
\begin{equation}
  \label{eq:orth-basis-decomp}
  W_m(g)=\scalar{V_m}{f_g}\quad \textrm{ where } g\in\MoebS \textrm{
    and }   m\in \spZ.
\end{equation} The explicit expression of \(W_m(g)\) could be derived
from the decomposition of the Cauchy kernel
in~\cite[\S~11.4.2]{BraDelSom82}, but it is important for us now that
formula~\eqref{eq:orth-basis-decomp} for a fixed \(g\) is a sort of
\emph{wavelet transform} \(\FSpace{H}{2}(\Space{S}{n-1})\rightarrow
\FSpace{C}{}(\spZ)\),
cf.~\eqref{eq:wavelet-transform-abstract}. We also use the
following properties of functions \(V_m(\vecbf{x})\) related to the
representation theory:
\begin{enumerate}
\item Functions \(V_m(\vecbf{x})\) with fixed
  \(\modulus{m}=m_1+\cdots+m_n\) form an \(O(n)\)-invariant irreducible
  module~\cite[\S~3.3]{GilbertMurray91}, which is required by the
  general construction of Taylor series~\cite[\S~3.4]{Kisil97c}.
\item There is the set of creation \(a^{+}_j\) and annihilation
      \(a^{-}_j\) operators (known
  from quantum mechanics):
  \begin{eqnarray*}
    a^{+}_j:& V_m(\vecbf{x})\mapsto V_{m'}(\vecbf{x}),&\qquad \textrm{
      where }
    m'=(m_1,\ldots, m_j+1, \ldots,m_n).\\
    a^{-}_j:& V_m(\vecbf{x})\mapsto m_jV_{m'}(\vecbf{x}),&\qquad
    \textrm{ where }
    m'=(m_1,\ldots, m_j- 1, \ldots,m_n).
  \end{eqnarray*}
  These operators satisfied~\cite{CnopsKisil97a} to the Heisenberg
  commutation relations:
  \begin{displaymath}
    [a^+_j,a^-_k]=\delta_{j,k}I, \qquad [a^+_j,a^+_k]=0, \qquad
    [a^-_j,a^-_k]=0. 
  \end{displaymath} Thus we have~\cite{CnopsKisil97a} a representation
  of the Heisenberg group \(\Space{H}{n}\) in
  \(\FSpace{H}{2}(\Space{S}{n-1})\).  Note also that
  in~\cite{KisRam96a} operators \(A^+_j\) and \(A^-_j\)were associated
  with operator of ``\(\times\)-product'' with the hypercomplex
  variable \(e_1x_j-e_jx_1\) partial derivative \(\partial_j\)
  correspondingly. 
\item The function \(V_0(\vecbf{x})\equiv 1\) coincides with the
  vacuum vector \(f_0(\vecbf{x})\) and:
  \begin{displaymath}
    V_m(\vecbf{x})=(a^+_1)^{m_1}(a^+_2)^{m_2}\cdots(a^+_n)^{m_n}
    f_0(\vecbf{x}). 
  \end{displaymath}
\end{enumerate}
Clearly we can decompose any shifted function
\(\uir{1}(g)V_k(\vecbf{x})\) over the basis \(V_m(\vecbf{x})\) in a
way similar to~\eqref{eq:Cauchy-decomposition}:
\begin{equation}
  \label{eq:rho-through-Zm}
  \uir{1}(g)V_k(\vecbf{x})=\sum_{m\in\spZ}
  W_{k,m}(g)\,V_m(\vecbf{x}), 
  \quad \textrm{ where } W_{k,m}(g)=\scalar{V_m}{\uir{1}(g)V_k}.
\end{equation} The representation
property \(\uir{1}(g_1)\uir{1}(g_2)=\uir{1}(g_1g_2)\) implies an
addition formula:
\begin{equation}
  \label{eq:token-property}
  W_{l,m}(g_1g_2)=\sum_{l\in\spZ} W_{l,k}(g_1)\,
  W_{k,m}(g_2). 
\end{equation} Thus functions \(W_{k,m}(g)\) are
\emph{tokens}~\cite{Kisil97b,Kisil01b} from the cancellative semigroup
\(\spZ\) to \(\MoebS\). This means that the
formula~\eqref{eq:rho-through-Zm} defines the representation
\(\uir{1}\) of \(\MoebS\) through the convolution on
\(\spZ\).

\section{Representations of $\MoebS$ in Algebras and Moduli}

A simple but important observation is that the M\"obius
transformations~\eqref{eq:moebius} can be easily extended to some
non-commutative \(\Cstar\)-algebras.

Let \(\algebra{A}\) be a \(\Cstar\)-algebra with the unit \(I\),
and an \(n\)-tuple \(A\) of  self-adjoint elements \(A_j\in\algebra{A}\),
\(j=1,\ldots,n\) be fixed. 
We consider the tensor product \(\algebra{A}\otimes\Cliff{n}\), which
we \emph{denote by} \(\algebra{A}_n\) for the brevity.  Its unit
element will be again denoted by \(I\). Then the tuple \(A\) can be
associated with the element \( \vecbf{A}=e_1A_1+e_2A_2+\cdots +e_nA_n
\) in \(\algebra{A}_n\). Let \(M\) be a left normed
\(\algebra{A}\)-module,  We denote by \(M_n\) the tensor product
\(M\otimes\Cliff{n}\). \(M_n\)  is a left \(\algebra{A}_n\) module of
course.  All constructed functional calculi according
to Definition~\ref{de:functional-calculus-new} are \(M_n\)-valued.

\subsection{Resolvent Approach}
\label{sec:resolvent-approach}

We define an action of the M\"obius group \(\Moeb{n}\) on the algebra
\(\algebra{A}_n\) by the natural formula (in
Notation~\ref{no:division}) similarly to expression~\eqref{eq:sp-rep}: 
\begin{equation}
  \label{eq:moebius-on-A}
  g: \vecbf{A} \mapsto g^{-1}\vecbf{A}=\frac{\bar{a}  \vecbf{A}
    -\bar{b} I}{a^*\vecbf{A}-b^* I}, 
  \qquad
  \phi_{(u,w)}=    \begin{pmatrix}
    \bar{a} & -\bar{b} \\ -b^*& a^*
  \end{pmatrix}
  \in \Moeb{\Space{S}{n-1}}.
\end{equation} To this end we need invertibility of the operator
\(a^*\vecbf{A} -b^*I\) in \(\algebra{A}_n\), which due to
invertibility of \(a^*\) in \(\Cliff{n}\) is equivalent to
invertibility of \(\vecbf{A} -\vecbf{u} I\), where \(\vecbf{u}
=(a^*)^{-1}b^*=a'b^*/\modulus{a}^2\) and thus
\(\modulus{\vecbf{u}}<1\). Therefore we arrive to the following
definition:
\begin{defn} \cite[Defn.~3.1]{Kisil95i}
  \label{de:Cliff-alg-spectr}
  The \emph{Clifford (algebraic) resolvent set} \(R( \vecbf{A})\) of an
  \(n\)-tuple \(A_1\), \(A_2\), \ldots, \(A_n\) is the maximal open
  subset of \(\Space{R}{n}\) such that for
  \(\vecbf{u}=u_1e_1+u_2 e_2+\cdots+u_n e_n\in
  R( \vecbf{A})\) 
  the element \(\vecbf{A}-{u} I\) is invertible in
  \(\algebra{A}_n\). 

  The \emph{Clifford (algebraic) spectrum} is the completion of the
  Clifford resolvent set \(\Space{R}{n}\setminus R({A})\).
\end{defn} 
\begin{rem}
  The Clifford (algebraic) spectrum is mainly an abbreviation for ``the
  complement of the resolvent set'' rather than an important
  characterisation of operator \(\vecbf{A}\). Such a characterisation
  is provided instead by the spectrum, defined through the support of
  functional calculus, see below. 
\end{rem}
Under the assumption that the Clifford algebraic spectrum of
\(\vecbf{A}\) belongs to the open unit ball \(\Space{B}{n}\) the orbit
\(\Space{A}{}=\{g^{-1} \vecbf{A} \such g\in \MoebS\}\) is a well
defined subset of \(\algebra{A}_n\). As any orbit \(\Space{A}{}\) is
a \(\MoebS\)-homogeneous space.

\begin{lem}[\textrm{\cite[Lem.~3.18]{Kisil95i}}]
  \label{le:mobius-tr-spec}
  For \(g\in\MoebS\) such that \((a^*)^{-1}b^*\in R(
  \vecbf{A})\) we have: 
  \begin{displaymath}
    \frac{\bar{a}\vecbf{A}-\bar{b}I}{a^*I-
      b^*\vecbf{A}}-\frac{\bar{a}\vecbf{x}-\bar{b}I}{a^*-b^*\vecbf{x}}= 
    (a-\vecbf{x}^*b)^{-1} (\vecbf{A}-\vecbf{x}I)(a^*I-b^*\vecbf{A})^{-1} .
  \end{displaymath}
  Consequently 
  \(    \vecbf{x}\in R(  \vecbf{A})\) implies 
    \(\frac{\bar{a}\vecbf{x}-\bar{b}}{a^*-b^*\vecbf{x}}\in
    R\left(\frac{\bar{a}\vecbf{A}-\bar{b}I}{a^*I- b^*\vecbf{A}}\right)\).
\end{lem}
\begin{proof}
  M\"obius transforms of vectors are vectors, for them
  \(\vecbf{y}^*=\vecbf{y}\), thus we have:
  \begin{eqnarray*}
    \label{eq:moeb-tr-spec}
    \lefteqn{\frac{\bar{a}\vecbf{A}-\bar{b}I}{a^*I-
      b^*\vecbf{A}}-\frac{\bar{a}\vecbf{x}-\bar{b}}{a^*-b^*\vecbf{x}}
    = \frac{\bar{a}\vecbf{A}-\bar{b}I}{a^*I-
      b^*\vecbf{A}}-\left(\frac{\bar{a}\vecbf{x}-\bar{b}}{a^*-b^*\vecbf{x}}\right)^*} 
  \qquad\\
    \qquad&=& (\bar{a}\vecbf{A}-\bar{b}I)(a^*I- b^*\vecbf{A})^{-1}-
    (a-\vecbf{x}b)^{-1}(\vecbf{x}a'-{b'})\\
    &=& (a-\vecbf{x}b)^{-1}\left((a-\vecbf{x}b)(\bar{a}\vecbf{A}-\bar{b}I)-
    (\vecbf{x}a'-{b'})(a^*I- b^*\vecbf{A})\right)(a^*I- b^*\vecbf{A})^{-1}\\
    &=&(a-\vecbf{x}^*b)^{-1} (\vecbf{A}-\vecbf{x}I)(a^*I-b^*\vecbf{A})^{-1} .
  \end{eqnarray*}
  The second statement follows from that result immediately.
\end{proof}

We define the \emph{resolvent} function
\(R(g,\vecbf{A}):\MoebS \times \Space{A}{}\rightarrow \algebra{A}_n\) by the
familiar expression: 
\begin{displaymath}
  R(g, \vecbf{A})=(a^* I-b^* \vecbf{A})^{-1} \quad 
\end{displaymath}
then a direct calculation shows that
\begin{equation}
  \label{eq:ind-rep-multipl}
  R(g_1,\ga)R(g_2,g_1^{-1}\ga)=R(g_1g_2,\ga).
\end{equation} The last identity is well known in representation
theory~\cite[\S~13.2(10)]{Kirillov76} and is a key ingredient of
\emph{induced representations}. Thus we can again
linearise~\eqref{eq:moebius-on-A} (cf.~\eqref{eq:rho-1-1}) in a
suitable space of \(M_n\) valued functions, where \(M_n\) is a left
\(\algebra{A}_n\) module as discussed at the beginning of this
section. We linearise~\eqref{eq:moebius-on-A} in the space of
continuous functions \(\FSpace{C}{}(\Space{A}{},M_n)\) as follows:
\begin{eqnarray}
  \uir{\vecbf{A}}(g_1): f(g^{-1} \vecbf{A} ) &\mapsto&
  R(g_1^{-1}g^{-1},  \vecbf{A})\,f(g_1^{-1}g^{-1}  \vecbf{A}) \label{eq:rho-a}\\
  &&\quad =
  (a^* I-b^* \vecbf{A})^{-1} \,
  f\!\left(
    \frac{\bar{a} \vecbf{A} - \bar{B}I}{a^*  I - b^*  \vecbf{A}} 
  \right).  \nonumber
\end{eqnarray} However such a representation is not unitary in
\(\FSpace{H}{2}(\MoebS)\) for \(n>2\) as can be seen from a
comparison with~\eqref{eq:rho-1-1}. To fix this we need an operator
which is symbolically represented by \(\modulus{a^*I-b^*\vecbf{A}}^{-n}\). 
When all operators \(A_j\) commute each other we can simply define:
\begin{eqnarray}
\modulus{a^*I-b^*\vecbf{A}}^{-2}&=&(a^*I-b^*\vecbf{A})^{-1}(a'I-\vecbf{A}b')^{-1}= 
\left(\modulus{a}^2+\modulus{b}^2\vecbf{A}^{2}\right)^{-1}\nonumber \\
&=&
\bigg(\modulus{a}^2-\modulus{b}^2\sum_{j=1}^n A_j^2\bigg)^{-1}.
\label{eq:modulus-of-op}
\end{eqnarray}
Then for an even \(n\geq 4\) we can straightforwardly define
\(\modulus{a^*I-b^*\vecbf{A}}^{-n+2}\). For an odd \(n\) we can define
a square root of the selfadjoint element
\((a'I-\vecbf{A}b')(a^*I-b^*\vecbf{A})\) of \(\algebra{A}_n\) by
various means. However the Clifford algebraic spectrum may not
guarantee the invertibility in~\eqref{eq:modulus-of-op}, thus some
additional assumptions of the type
\(\norm{\vecbf{A}}<(1+\sqrt{2})^{-1}\) are
required~\cite{JefferiesMcIntosh99a}. 
Consequently 
for  for a commuting \(n\)-tuple
(\(n>2\)) of operators \(A_j\) one defines a
representation \(\uir{\vecbf{A}}\) in the
\(\FSpace{C}{}(\Space{A}{},M_n)\) by the expression:
\begin{equation}
  \label{eq:repr-mcIntosh}
  \uir{\vecbf{A}} f(\vecbf{A}) = R(g,\vecbf{A})
  \modulus{a^*I-b^*\vecbf{A}}^{-n+2} f(g^{-1}\vecbf{A}).
\end{equation}
In this way we obtain the monogenic calculus of commuting operators
studied in~\cite{McInPryde87}. 

For any \(v\in M\) we can again define a \(K\)-invariant
\emph{vacuum vector} as \(f_0(\vecbf{A},v)=v\otimes f_0(\vecbf{A})\equiv v \in
\FSpace{C}{}(\Space{A}{},M_n)\).  
It generates the associated with \(f_{\vecbf{u}}\) family of \emph{coherent
  states} \(f_g(\vecbf{A},v)=R(g,\vecbf{A})
  \modulus{a^*I-b^*\vecbf{A}}^{-n+2} v\),
where \(g\in\MoebS\). 
The \emph{wavelet transform}  defined by
the same common formula based on coherent states (cf.~\eqref{eq:cauchy}):
\[\oper{W}_m f(g)= \scalar{\uir{\vecbf{A}}(g) f_0}{f},\]
is a version of Cauchy integral, which maps
\(\FSpace{L}{2}(\Space{A}{},M')\) to
\(\FSpace{C}{}(\MoebS,\Space{C}{})\), where \(M'\) is the dual of the
module \(M\). The classical Riesz-Dunford functional calculus is a
particular realisation of this approach~\cite{Kisil02a}.

For a non-commuting tuple \(\vecbf{A}\) one can, for example, define
representation \(\uir{\vecbf{A}}\) using the fruitful
approach~\cite{JefferiesMcIntosh99a} based on the plain wave
decomposition~\cite{Sommen88}.  An alternative is the Taylor expansion
construction initiated in~\cite{Kisil95i}.

\subsection{Taylor Expansion Approach}
\label{sec:tayl-expans-appr}
To define a functional calculus for \(\vecbf{A}\) we fix images
\(\Phi(V_m)\) of \(V_m\) (see 
Subsection~\ref{sec:taylor-series}), \(m\in\spZ\) in
\(\algebra{A}_n\), cf.~\cite{KisRam95a,KisRam96a}. Seemingly this
could be done in many different 
ways, but the covariance property fixes one preferred assignment. 
Indeed subgroup \(O(n)\) of \(\MoebS\) contains permutations of
elements of orthonormal basis \(e_k\).
Functions \(V_m(\vecbf{x})\)
are symmetric polynomials of \(x_j\) and are invariant under such
permutations. To preserve \(O(n)\) invariance we define
\(\Phi(V_m)=\Phi_{\vecbf{A},x}(V_m)\) associated to the tuple
\(\vecbf{A}\)  to be 
\begin{equation}
  \label{eq:Am-defn}
  \Phi(V_m)=A_m :=\frac{1}{\modulus{m}!}\sum_{\sigma\in
    S_{\modulus{m}}} e_{\sigma(1)}A_{\sigma(1)} 
  e_{\sigma(2)}A_{\sigma(2)}\cdots
  e_{\sigma(n)}A_{\sigma(\modulus{m})} 
\end{equation}
is the averaging of products of \(m_j\) copies of \(e_jA_j\) 
 over the permutation group \(S_{\modulus{m}}\). 

The  value \(r_R(\vecbf{A})=\lim_{j\rightarrow\infty}\sup_\sigma 
\norm{A_{\sigma(1)} \cdots A_{\sigma(j)} }^{1/j}\), \(1\leq\sigma(i)\leq
n\) is known as the \emph{Rota-Strang joint spectral
radius}~\cite{Rota60a}. We give a similar definition which is better
tailored to our circumstances:
\begin{defn}
  Let \(m\in\spZ\), \(v\in M\) and \(A_m\) be defined
  in~\eqref{eq:Am-defn}.  We call 
  \begin{displaymath}
    r_S(\vecbf{A})=\limsup_{\modulus{m}\rightarrow\infty}
    \norm{A_m}^{1/\modulus{m}}_\algebra{A} \qquad \textrm{ and } \qquad
    r_L(\vecbf{A},v)=\limsup_{\modulus{m}\rightarrow\infty}
    \norm{A_m v}^{1/\modulus{m}}_M
  \end{displaymath} the \emph{symmetric joint spectral
    radius} of \(\vecbf{A}\) and \emph{local spectral radius} of \(\vecbf{A}\) at \(v\)
  correspondingly. Obviously \(r_S(\vecbf{A})\leq r_R(A)\) and
  \(r_L(\vecbf{A},v)\leq r_S(A)\norm{v}\).
\end{defn}
Let \(r_L(\vecbf{A},v)< 1\), \(v\in M\) and a sequence \(c_m\),
\(m\in\spZ\) be a square summable. Then the infinite
series \(\sum_{m\in\spZ} c_m A_m v\) is absolutely convergent by norm
in  \(M_n\). The linear space of all such sequences is denoted by
\(\FSpace{H}{2}(\vecbf{A},v)\).
Analogously to representation \(\uir{1}\)
in~\eqref{eq:rho-through-Zm} we define
 an action \(\uir{\vecbf{A},v}\) of \(\MoebS\) on \(\FSpace{H}{2}(A,v)\) by:
 \begin{equation}
   \label{eq:rho-A-defn-taylor}
   \uir{\vecbf{A},v}(g): \sum_{k\in\spZ} c_k A_k v \mapsto \sum_{k\in\spZ} d_k
   A_k v, \quad \textrm{ where }
   d_k=\sum_{m\in\spZ}
  W_{k,m}(g)\,c_m.
 \end{equation}
 Then the identity~\eqref{eq:token-property} implies that
 \(\uir{\vecbf{A},v}\) is a representation of \(\MoebS\). 

\begin{defn}
  \label{de:monogenic-calculus}
  Let \(r_L(A,v)< 1\) then the monogenic functional calculus
  \(\Phi=\Phi_{\vecbf{A},v}\) associated to a \(n\)-tuple \(\vecbf{A}\) and a vector
  \(v\in M\) is a continuous linear map \(\Phi:
  \FSpace{H}{2}(\Space{S}{n-1}) \rightarrow \FSpace{H}{2}(A,v)\) is
  defined by the following two conditions:
  \begin{enumerate}
  \item \label{item:intertwining}
    \(\Phi\) intertwines \(\uir{1}\)~\eqref{eq:rho-1-1} and
    \(\uir{\vecbf{A},v}\)~\eqref{eq:rho-A-defn-taylor}: \(\Phi\uir{1}(g)=
    \uir{\vecbf{A},v}(g)\Phi\) for all 
    \(g\in\MoebS\).  
  \item 
    \label{item:initial}
    The map of  vacuum vectors is \(\Phi(f_0)=f_v\), where
    \(f_0(\vecbf{x})\equiv 1\) and \(f_v= v\).
  \end{enumerate}
\end{defn}
This defines monogenic calculus uniquely,
particularly its integral formula.
\begin{prop}
  Let \(E(g,\vecbf{A})\) be the family of coherent states for \(\uir{\vecbf{A},v}\):
  \begin{eqnarray}
    \label{eq:cauchy-kernel-operator}
    E(g,\vecbf{A})=\uir{\vecbf{A},v} f_v= \sum_{k\in\spZ} W_{m,0}(g) A_m v.
  \end{eqnarray}
  Then the functional calculus \(\Phi_{\vecbf{A},v}\) is defined by the
  Integral formula:
  \begin{displaymath}
    \Phi_{\vecbf{A},v}f = \int_{\MoebS} E(g,\vecbf{A}) f(g)\,dg.
  \end{displaymath}
\end{prop}
\begin{proof}
  Indeed using the Definition~\ref{de:monogenic-calculus} we
  calculate for \(f=\scalar{\uir{1}(g) f_0}{f}\):
  \begin{eqnarray}
    \Phi_{\vecbf{A},v}f &=& \Phi_{\vecbf{A},v} \scalar{\uir{1}(g) f_0}{f} \nonumber \label{eq:calc-tr1}\\
   &=&  \scalar{\Phi_{\vecbf{A},v} \uir{1}(g) f_0}{f} \label{eq:calc-tr2}\\
   &=&  \scalar{\uir{\vecbf{A},v}(g) \Phi_{\vecbf{A},v}  f_0}{f} \label{eq:calc-tr3}\\
   &=&  \scalar{\uir{\vecbf{A},v}(g) f_v}{f} \label{eq:calc-tr4}\\
   &=&  \scalar{E(g,\vecbf{A})}{f} \label{eq:calc-tr5}\\
   &=& \int_{\Space{S}{n-1}} E(g,\vecbf{A}) f(\vecbf{x})\,d\vecbf{x},\nonumber
  \end{eqnarray}
  where~\eqref{eq:calc-tr2} is obtained by linearity and continuity of
  functional calculus, \eqref{eq:calc-tr3} follows from the
  intertwining property~\ref{item:intertwining}, \eqref{eq:calc-tr4}
  is obtained from the initialisation property~\ref{item:initial}, and
  finally~\eqref{eq:calc-tr5} uses
  expression~\eqref{eq:cauchy-kernel-operator} for \(E(g,\vecbf{A})\).
\end{proof}

The full consideration of the monogenic calculus and the corresponding
joint spectrum requires a solid background from the representation
theory of semisimple Lie groups. We will consider a simpler but still
illustrative case in the next section.

\section{Functional Calculus and Spectrum for a Pair of Matrices}
\label{sec:funct-calc-spectr}

In this section we demonstrate the previous construction by the
simplest non-trivial example: functional calculus for a pair \(A_1\),
\(A_2\) of self-adjoint non-commuting operators with finite
dimensional ranges, cf.~\cite{JefferiesStraub03}. Instead of tensor product
\(e_1A_1+e_2A_2\) with Clifford algebra \(\Cliff{2}\) we can consider
the complexification \(A_1+iA_2\) since the product \(i=e_1e_2\) has
all properties of the complex imaginary unit. The group \(\MoebS\) is
the \(\SL\) group in this case, \(O(2)\) consists from the orthogonal
rotations of the plane, and \(\MoebS/O(n)=\SL/O(2)\) is the unit
disk \(\Space{D}{}\). 
In two dimensions the formula~\eqref{eq:rho-a} defines a isometric
representation of \(\MoebS\) without a normalising factor. 

\subsection{Jet Bundles and Prolongations of $\uir{1}$}
\label{sec:jet-bundl-prol-1}
To formulate the complete description of monogenic calculus and
spectrum of \(\vecbf{A}=A_1+iA_2\) we use the language of jet
spaces and prolongations of representations introduced by S.~Lie, see
\cite{Olver93,Olver95} for a detailed exposition.

\begin{defn} \textup{\cite[Chap.~4]{Olver95}}
  Two holomorphic functions have \(n\)th \emph{order contact} in a point
  if their value and their first \(n\) derivatives agree at that point,
  in other words their Taylor expansions are the same in first \(n+1\)
  terms. 

  A point \((z,u^{(n)})=(z,u,\vecbf{u}_1,\ldots,\vecbf{u}_n)\) of the \emph{jet space}
  \(\Space{J}{n}\sim\Space{D}{}\times\Space{C}{n}\) is the equivalence
  class of holomorphic functions having \(n\)th contact at the point \(z\)
  with the polynomial:
  \begin{equation}\label{eq:Taylor-polynom}
    p_n(w)=\vecbf{u}_n\frac{(w-z)^n}{n!}+\cdots+\vecbf{u}_1\frac{(w-z)}{1!}+u.
  \end{equation}
\end{defn}

For a fixed \(n\) each holomorphic function
\(f:\Space{D}{}\rightarrow\Space{C}{}\) has \(n\)th \emph{prolongation}
(or \emph{\(n\)-jet}) \(\object[_n]{j}f: \Space{D}{} \rightarrow
\Space{C}{n+1}\): 
\begin{equation}\label{eq:n-jet}
  \object[_n]{j}f(z)=(f(z),f'(z),\ldots,f^{(n)}(z)).
\end{equation}The graph \(\Gamma^{(n)}_f\) of \(\object[_n]{j}f\) is a
submanifold of \(\Space{J}{n}\) which is section of the \emph{jet
bundle} over \(\Space{D}{}\) with a fibre \(\Space{C}{n+1}\). We also
introduce a notation \(J_n\) for the map \(
  J_n:f\mapsto\Gamma^{(n)}_f
\) of a holomorphic \(f\) to the graph \(\Gamma^{(n)}_f\) of its \(n\)-jet
\(\object[_n]{j}f(z)\)~\eqref{eq:n-jet}.

One can prolong any map of functions \(\psi: f(z)\mapsto [\psi f](z)\) to
a map \(\psi^{(n)}\) of \(n\)-jets by the formula
\begin{equation}\label{eq:prolong-def}
  \psi^{(n)} (J_n f) = J_n(\psi f).
\end{equation} For example such a prolongation \(\uir{1}^{(n)}\) of the
representation \(\uir{1}\) of the group \(\MoebS\) in
\(\FSpace{H}{2}(\Space{D}{})\) (as any other representation of a Lie
group~\cite{Olver95}) will be again a representation of
\(\MoebS\). Equivalently we can say that \(J_n\) \emph{intertwines} \(\uir{1}\) and
\(\uir{1}^{(n)}\):
\begin{displaymath}
   J_n \uir{1}(g)= \uir{1}^{(n)}(g) J_n \quad
  \textrm{ for all } g\in\MoebS.
\end{displaymath}
Of course, the representation \(\uir{1}^{(n)}\) is not irreducible: any jet
subspace \(\Space{J}{k}\), \(0\leq k \leq n\) is
\(\uir{1}^{(n)}\)-invariant subspace of \(\Space{J}{n}\).  However the
representations \(\uir{1}^{(n)}\) are
\emph{primary}~\cite[\S~8.3]{Kirillov76} in the sense that they are not 
sums of two subrepresentations.

The following statement explains why jet spaces appeared in our study.
\begin{prop}
  \label{pr:Jordan-zero} Let the matrix \(\vecbf{A}=A_1+iA_2\) be a
  Jordan block of a length \(k\) with the eigenvalue \(u=0\), and
  \(v\) be its root vector of order \(k\), i.e. \(\vecbf{A}^{k-1}v\neq
  \vecbf{A}^k v =0\). Then the restriction of \(\uir{\vecbf{A},v}\) on the
  subspace generated by \(v_m\) is equivalent to the representation
  \(\uir{1}^{k}\).
\end{prop}

\subsection{Spectrum and the Jordan Normal Form of a Matrix} Now we
are prepared to describe a spectrum of a matrix
\(\vecbf{A}=A_1+iA_2\). Since the functional calculus is an
intertwining operator its support is a decomposition into intertwining
operators with prime representations (we could not expect generally
that these prime subrepresentations are irreducible).

Recall the group of inner automorphisms \(b_g: g_1\mapsto b_g
(g_1)=g^{-1}g_1g\) of \(\MoebS\). The representation
\(\uir{g}(g_1)=\uir{1}(b_g(g_1))\) is equivalent to \(\uir{1}\) and
they are obviously intertwined by the operator \(\uir{1}(g^{-1})\):
\(\uir{g}\uir{1}(g^{-1})=\uir{1}(g^{-1})\uir{1}\). For a Jordan block
\(\vecbf{A}\) with an eigenvalue \(\vecbf{u}\) its M\"obius
transformation with the matrix \(\matr{1}{\vecbf{u}'}{\vecbf{u}}{1}\)
will be a Jordan block with eigenvalue \(\vecbf{0}\) due to
Lemma~\ref{le:mobius-tr-spec}. Thus inner automorphisms  extend
Proposition~\ref{pr:Jordan-zero} to the complete 
characterisation of \(\uir{\vecbf{A},v}\) for matrices.
\begin{prop}  
  \label{pr:3d-spectr}
  Representation \(\uir{\vecbf{A},v}\) is equivalent to a direct sum of the
  prolongations \(\uir{1}^{(k)}\) of \(\uir{1}\) in the \(k\)th jet space
  \(\Space{J}{k}\) intertwined with inner automorphisms. Consequently
  the \textit{spectrum} of \(\vecbf{A}\) (defined via the functional calculus
  \(\Phi_{\vecbf{A},v}\)) labelled exactly by \(n\) pairs of numbers
  \((\vecbf{u}_i,k_i)\), where \(\vecbf{u}_i\in\Space{D}{}\),
  \(k_i\in\Space[+]{Z}{}\) for \(1\leq i \leq n\) some of whom could
  coincide.
\end{prop}
Obviously this spectral theory is a fancy restatement of the \emph{Jordan
  normal form} of matrices.

\begin{figure}[tb]
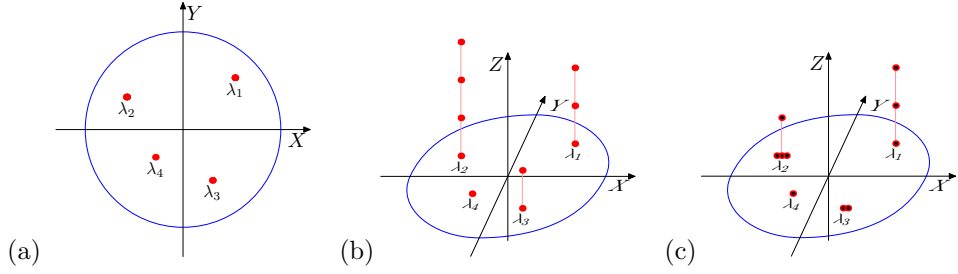

  \begin{center}
 (a) \includegraphics[scale=.65]{calc1vr.10}\hfill
  (b)\includegraphics[scale=.65]{calc1vr.27}\hfill
  (c)\includegraphics[scale=.65]{calc1vr.40}
 \caption[Three dimensional spectrum]{Classical spectrum (a) of a
   matrix vs. the new version (b) with its mapping (c).}
    \label{fig:3dspectrum}
  \end{center}
\end{figure}

\begin{example}
  \label{ex:3dspectrum}
  Let \(J_k(\vecbf{u})\) denote the Jordan block of the length \(k\) for the
  eigenvalue \(\vecbf{u}\). On the Fig.~\ref{fig:3dspectrum} there are two
  pictures of the spectrum for the matrix
  \begin{displaymath}
    a=J_3\left(\vecbf{u}_1\right)\oplus     J_4\left(\vecbf{u}_2\right) 
    \oplus J_1\left(\vecbf{u}_3\right) \oplus      J_2\left(\vecbf{u}_4\right),
  \end{displaymath} 
  where
  \begin{displaymath}
    \vecbf{u}_1=\frac{3}{4}e^{i\pi/4}, \quad
    \vecbf{u}_2=\frac{2}{3}e^{i5\pi/6}, \quad
    \vecbf{u}_3=\frac{2}{5}e^{-i3\pi/4}, \quad
    \vecbf{u}_4=\frac{3}{5}e^{-i\pi/3}.
  \end{displaymath} Part (a) represents the conventional two-dimensional
  image of the spectrum, i.e. eigenvalues of \(\vecbf{A}\), and
  \href{http://maths.leeds.ac.uk/~kisilv/calc1vr.gif}{(b) describes
  spectrum \(\spec{} a\)} arising from the wavelet construction. The
  first image does not allow to distinguish \(\vecbf{A}\) from many other
  essentially different matrices, e.g. the diagonal matrix
  \begin{displaymath}
    \diag\left(\vecbf{u}_1,\vecbf{u}_2,\vecbf{u}_3,\vecbf{u}_4\right),
  \end{displaymath}
  which even have a different dimensionality.
  At the same time the Fig.~\ref{fig:3dspectrum}(b)
  completely characterise \(\vecbf{A}\) up to a similarity. Note that each point of
  \(\spec \vecbf{A}\) on Fig.~\ref{fig:3dspectrum}(b) corresponds to a particular
  root vector, which spans a primary subrepresentation.
\end{example}
In light of the previous discussions~\cite[p.~29]{Kisil95i},
\cite[Ex.~6.3]{JeffMcInt98a} the following simple example is still of interest.  
\begin{example}
  For a pair Pauli matrices \(J_1=\matr{1}{0}{-1}{0}\) and
  \(J_2=\matr{0}{1}{1}{0}\), the joint Clifford algebraic spectrum
  (Definition~\ref{de:Cliff-alg-spectr}) as found
  in~\cite[p.~29]{Kisil95i} is the single point
  \((0,0)\in\Space{D}{}\). The joint spectrum found in
  \cite[Ex.~6.3]{JeffMcInt98a} coincides with the Weyl joint spectrum
  and the numerical range~\cite{JefferiesStraub03}: all of them are
  the entire unit disk \(\Space{D}{}\). Finally the joint spectrum
  from Proposition~\ref{pr:3d-spectr} is a pair of points
  \((\vecbf{0},0)\) and  \((\vecbf{0},1)\) from
  \(\Space{R}{2}\times\Space[+]{Z}{}\) since \(J_1+iJ_2\) is similar
  to the  Jordan block of the length \(2\) with the eigenvalue \(0\). 
\end{example}

\subsection{Spectral Mapping Theorem} As was mentioned in the
Introduction a reasonable spectrum should be linked to the
corresponding functional calculus by an appropriate spectral mapping
theorem. The new version of spectrum is based on prolongation of
\(\uir{1}\) into jet spaces. Naturally a correct version of spectral
mapping theorem should operate in jet spaces as well. 

Let \(\phi: \Space{D}{} \rightarrow \Space{D}{}\) be a holomorphic
map, let us define its action on functions \([\phi_*
f](\vecbf{z})=f(\phi(\vecbf{z}))\). According to the general
formula~\eqref{eq:prolong-def} we can define the prolongation
\(\phi_*^{(n)}\) onto the jet space \(\Space{J}{n}\). Its associated
action \(\uir{1}^k \phi_*^{(n)}=\phi_*^{(n)}\uir{1}^n\) on the pairs
\((\vecbf{u},k)\) is given by the formula:
\begin{equation}
  \label{eq:phi-star-action}
  \phi_*^{(n)}(\vecbf{u},k)=\left(\phi(\vecbf{u}),
    \left[\frac{k}{\deg_\vecbf{u} \phi}\right]\right),
\end{equation}
where \(\deg_\vecbf{u} \phi\) denotes the degree of zero of the function
\(\phi(\vecbf{z})-\phi(\vecbf{u})\) at the point
\(\vecbf{z}=\vecbf{u}\) and \([x]\) denotes 
the integer part of \(x\). We are ready to state

\begin{thm}[Spectral mapping] 
  Let \(\phi\) be a holomorphic mapping  \(\phi: \Space{D}{}
  \rightarrow \Space{D}{}\) and its prolonged action \(\phi_*^{(n)}\) defined
  by~\eqref{eq:phi-star-action}, then
  \begin{displaymath}
    \spec \phi(\vecbf{A}) = \phi_*^{(n)} (\spec \vecbf{A}). 
  \end{displaymath}
\end{thm} 

The explicit expression of~\eqref{eq:phi-star-action} for
\(\phi_*^{(n)}\), which involves derivatives of \(\phi\) up to \(n\)th order,
is known, see for
example~\cite[Thm.~6.2.25]{HornJohnson94}. However it was
not recognised before as a form of spectral mapping.

\begin{example}
  Let us continue with Example~\ref{ex:3dspectrum}. Let \(\phi\) map
  all four eigenvalues \(\vecbf{u}_1\), \ldots, \(\vecbf{u}_4\) of the
  matrix \(\vecbf{A}\) into themselves. Then Fig.~\ref{fig:3dspectrum}(a) will
  represent the classical spectrum of \(\phi(a)\) as well as \(\vecbf{A}\).
  In the contrast Fig.~\ref{fig:3dspectrum}(c) shows mapping of the new
  spectrum for the case
  \(\phi\)  has
  \textit{orders of zeros} at these points as follows: the order \(1\)
  at \(\vecbf{u}_1\), exactly the order \(3\) at \(\vecbf{u}_2\), an order
  at least \(2\) at \(\vecbf{u}_3\), and finally any order at
  \(\vecbf{u}_4\).
\end{example}

\bibliographystyle{amsplain} 
\bibliography{abbrevmr,akisil,analyse,algebra,aphysics}

\end{document}